\documentclass[a4paper,10pt]{amsart}

\usepackage{amsmath,amssymb,amsthm,a4wide}

\usepackage{tikz}
\usetikzlibrary{shapes}
\usetikzlibrary{intersections}

\usepackage{graphicx}
\usepackage{graphics}

\newtheorem{theorem}{Theorem}

\newtheorem*{proposition}{Proposition}

\newcommand{\QQ}{\mathbb{Q}}

\begin{document}

\title[Rigidity for the Maximal Function]{A Rigidity Phenomenon for the \\Hardy-Littlewood Maximal Function}
\author{Stefan Steinerberger}
\address{Department of Mathematics, Yale University, 10 Hillhouse Avenue, New Haven, CT 06511, USA}
\email{stefan.steinerberger@yale.edu}

\begin{abstract} The Hardy-Littlewood maximal function $\mathcal{M}$ and the trigonometric function $\sin{x}$ are two central objects in harmonic analysis.
We prove that $\mathcal{M}$ characterizes $\sin{x}$ in the following way: let $f \in C^{\alpha}(\mathbb{R}, \mathbb{R})$ be a periodic function and $\alpha > 1/2$. If there exists
a real number $0 < \gamma < \infty$ such that the averaging operator
$$ (A_xf)(r) = \frac{1}{2r}\int_{x-r}^{x+r}{f(z)dz}$$
has a critical point in $r = \gamma$ for every $x \in \mathbb{R}$, then
$$f(x) = a+b\sin{(cx + d)} \qquad \mbox{for some}~a,b,c,d \in \mathbb{R}.$$ 
This statement can be used to derive a characterization of trigonometric functions as those nonconstant functions for which
the computation of the maximal function $\mathcal{M}$ is as simple as possible.
The proof uses the Lindemann-Weierstrass theorem from transcendental number theory.
\end{abstract}
\maketitle

\section{Introduction and main results}

\subsection{Introduction} Maximal functions are a central object in harmonic analysis; conversely, harmonic analysis is built up from trigonometric functions. We were motivated by the simple question whether a maximal function is able to 'recognize' a trigonometric function in any particular way. We focus on the centered Hardy-Littlewood 
maximal function on the real line 
$$(\mathcal{M}f)(x) = \sup_{r > 0}\frac{1}{2r}\int_{x-r}^{x+r}{|f(z)|dz}.$$
Classical results are the embedding $\mathcal{M}:L^1 \rightarrow L^{1, \infty}$, where the sharp constant is known \cite{mel}, as well as the embedding $\mathcal{M}:L^p \rightarrow L^p$ for
$1 < p \leq \infty$. The wealth of theory developed around maximal functions can no longer be succinctly summarized: we refer to the book of Stein \cite{stein} for the
classical theory and a survey of Wolff \cite{wo} on the Kakeya problem.

\subsection{The interval length function $r_f(x)$.} Usual questions around maximal functions are concerned with their size: since $(\mathcal{M}f)(x) \geq |f(x)|$ in all Lebesgue points, it
is of interest to understand mapping properties in $L^p$ spaces. Our question goes in an orthogonal direction: how complicated is the dynamical behavior of the 'maximal' intervals?
The question is not well-posed because there might be more than one interval centered at $x$ over which the average value of the function coincides with the maximal function: in these cases, we opt for taking the smallest such interval. Formally, we define the length function $r_f(x):\mathbb{R} \rightarrow \mathbb{R}_{\geq 0}$ associated to a periodic function $f$ by
$$r_f(x) = \inf \left\{r \geq 0: \frac{1}{2r}\int_{x-r}^{x+r}{f(z)dz} = \sup_{s > 0}{\frac{1}{2s}\int_{x-s}^{x+s}{f(z)dz}}\right\},$$
where the integral is to be understood as the point evaluation if $r = 0$. It is easy to see that $r_f$ is well-defined and finite for periodic functions.

\subsection{Main results.} The purpose of this paper is to study the situation, where for all $x \in \mathbb{R}$ both $r_f(x)$ and $r_{-f}(x)$ are either 0 or a fixed positive real number and to show that this characterizes the trigonometric
function. This theorem may be understood as a characterization of trigonometric functions by means of a dynamical aspect of the 
Hardy-Littlewood maximal function. It seems to have surprisingly little to do with traditional characterizations involving geometry, power series, differential
equations or spectral theory. Indeed, we failed to find a slick reduction to any of the classical characterizations and 
ended up needing tools from transcendental number theory.

\begin{theorem}  Let $f \in C^{\alpha}(\mathbb{R}, \mathbb{R})$ be a periodic function and $\alpha > 1/2$. There exists a positive number $\gamma > 0$ such that the averaging operator
$$ (A_xf)(r) = \frac{1}{2r}\int_{x-r}^{x+r}{f(z)dz}$$
has a critical point in $r = \gamma$ for every $x$ if and only if
$$f(x) = a+b\sin{(cx + d)} \qquad \mbox{for some}~a,b,c,d \in \mathbb{R}.$$ 
\end{theorem}
Theorem 1 is the strongest statement in this paper; it is relatively easy to deduce the following
statement, which formulates everything in terms of the complexity of the maximal intervals.
\begin{theorem}  Let $f \in C^{\alpha}(\mathbb{R}, \mathbb{R})$ be a periodic function and $\alpha > 1/2$. Then
$$ \left| \bigcup_{x \in \mathbb{R}}{\left\{r_f(x), r_{-f}(x)\right\}} \right| \leq 2$$
if and only if 
$$f(x) = a+b\sin{(cx + d)} \qquad \mbox{for some}~a,b,c,d \in \mathbb{R}.$$ 
\end{theorem}
We emphasize that we do not even know whether the statement remains true if the constant 2 is replaced by a larger positive integer (but conjecture
that it does).
Another way of stating Theorem 2 is as follows: suppose
the periodic function $f(x)$ does not change sign and that both $\mathcal{M}f$ and $\mathcal{M}(-f)$ can be computed by checking the average over an
interval of fixed interval and comparing it with point evaluation, i.e. suppose there exists a fixed number $0 < \gamma < \infty$
such that $f$ satisfies the equation
$$ (\mathcal{M}f)(x) = \max\left(|f(x)|, \frac{1}{2\gamma}\int_{x-\gamma}^{x+\gamma}{|f(z)|dz}\right) \qquad \mbox{for all}~x \in \mathbb{R}$$
and the same condition (with the same value $\gamma$) holds for $\mathcal{M}(-f)$, then 
$$f(x) = a+b\sin{(cx + d)} \qquad \mbox{for some}~a,b,c,d \in \mathbb{R}.$$ 

\subsection{A delay differential equation.} Perhaps the most natural first step after seeing Theorem 1 would be 
to try a combination of differentiation and algebraic manipulations to obtain an ordinary differential equation (with
the hope of it being $y'' + y =0$). As it turns out, this does not work and produces much more interesting results
instead. Differentiation in $r$ implies that
$$ 0 = \partial_r \frac{1}{2r}\int_{x-r}^{x+r}{f(z)dz}\big|_{r = \gamma} = -\frac{1}{2\gamma^2}\int_{x-\gamma}^{x+\gamma}{f(z)dz} + \frac{1}{2\gamma}(f(x+\gamma) + f(x-\gamma)).$$
Assuming $f \in C^1$, this equation can now be differentiated in $x$ and yields
$$ f'(x+\gamma) - \frac{1}{\gamma}f(x+\gamma) = -f'(x-\gamma) - \frac{1}{\gamma}f(x-\gamma).$$
The qualitative theory of delay differential equations is a lot more complicated than the theory of ordinary differential equations because the space of
solutions is \textit{much} larger (uncountable): any $C^1$ function $g:[0, 2\gamma] \rightarrow \mathbb{R}$ with correct boundary
conditions can always be extended to a solution of the delay differential equation. However, as an easy consequence of Theorem 2, we
can show that there are few periodic solutions.
\begin{theorem} Let $\gamma > 0$ be fixed and let $f \in C^1(\mathbb{R}, \mathbb{R})$ be a solution of the delay differential equation
$$ f'(x+\gamma) - \frac{1}{\gamma}f(x+\gamma) = -f'(x-\gamma) - \frac{1}{\gamma}f(x-\gamma).$$
If $f$ is periodic, then
$$f(x) = a+b\sin{(cx + d)} \qquad \mbox{for some}~a,b,c,d \in \mathbb{R}.$$
\end{theorem}
Considering the large (uncountable) number of solutions, it is utterly remarkable that there are so few periodic solutions. We have not been able to locate any type of argument in the literature 
that would allow to prove a result of this type.

\subsection{Open questions.} We believe that many of the assumptions can be weakened. Periodicity of the functions is necessary to allow the use of Fourier series on which
our argument is based, however, it seems natural to assume that the properties discussed could not hold for a nonperiodic function anyway. The assumption
$f \in C^{\alpha}(\mathbb{R}, \mathbb{R})$ with $\alpha > 1/2$ is required at one point in the proof to enforce uniform convergence of the Fourier series; again, it seems to
be an artefact of the method. We note that the condition $f \in C^{\alpha}(\mathbb{R}, \mathbb{R})$ with $\alpha > 1/2$ in our statements could always be replaced with the 
condition of $f$ having an absolutely convergent Fourier series. The strongest statement we believe could
to true is the following.
\begin{quote}
\textit{Conjecture.} If $f \in L^{\infty}(\mathbb{R})$ is real-valued and $r_f(x)$ assumes only finitely many different values, then 
 $$f(x) = a+b\sin{(cx + d)} \qquad \mbox{for some}~a,b,c,d \in \mathbb{R}.$$ 
\end{quote}
A more daring conjecture would be that it suffices to assume that
$$\bigcup_{x \in \mathbb{R}}{\left\{r_f(x)\right\}} \subset \mathbb{R}
 \qquad \mbox{is a Lebesgue-null set.}$$
 If $r_f(x)$ is contained in a set of 'small' non-zero Lebesgue measure, what does that imply for the function? It seems to indicate,
in some weak sense, that Fourier frequencies interact weakly (perhaps in the sense of a $\Lambda(p)-$property?). Furthermore, it
seems that if $f$ is given by a lacunary Fourier series, then
$$\bigcup_{x \in \mathbb{R}}{\left\{r_f(x)\right\}} \subset \mathbb{R} \qquad \mbox{can have 'fractal' structure}$$
again in a vague sense (small measure and a very large number of connected components): it could be of interest to try to understand quantitative versions of this
basic intuition. One could also ask to which extent this persists in higher dimensions: in $\mathbb{R}^d$, if we consider the maximal function associated to axis-parallel 
rectangles and the natural analogue $r_f(x):\mathbb{R}^d \rightarrow \mathbb{R}^d$, then setting
$$f(x_1, \dots, x_d) = \prod_{i=1}^{d}{(a_i + b_i\sin{x_i} + c_i\cos{x_i})} \quad \mbox{implies} \quad \left| \bigcup_{x \in \mathbb{R}}{\left\{r_f(x), r_{-f}(x)\right\}} \right| \leq 2^d.$$
Already in two dimensions, there are many natural maximal functions and
it is not clear to us whether similar statements hold for any of them. We recall the Pompeiu conjecture \cite{pomp}: if $K \subset \mathbb{R}^n$ is a simply connected Lipschitz domain and $f: \mathbb{R}^n \rightarrow \mathbb{R}$ is a nonzero
continuous function such that the integral of f vanishes over every congruent copy of $K$ -- does this imply that $K$ is a ball? Is there a connection between the Pompeiu conjecture
and the maximal problem for disks?\\

\textit{The discrete setting.} In light of the recent results \cite{bob, car, car1} concerning the behavior of the maximal function on the lattice, this might be another interesting direction to investigate. For a function $f:\mathbb{Z} \rightarrow \mathbb{R}$, we define the maximal function as
$$ (\mathcal{M}f)(n) = \sup_{r \in \mathbb{N}_{\geq 0}}{\frac{1}{2r +1}{\sum_{k=n-r}^{n+r}{f(k)}}}.$$
The length function $r_{f}:\mathbb{Z} \rightarrow \mathbb{N}$ is defined as above. Numerical experiments show that
the continuous case translates into the discrete setting: for generic parameters functions of the form 
$$ f(n) = a + b \sin{(c n + d)} \qquad a,b,c,d \in \mathbb{R}$$
seem to give rise to two-valued $r_f$. We do not have a formal proof of this statement, it should be equivalent to a series
of trigonometric inequalities that might actually be in the literature. The property is stable under small perturbations. However, there also exist completely
different functions with a two-valued $r_f$: the following example was found by Xiuyuan Cheng (personal communication). Taking
$$ f(n) = \frac{(-1)^n}{\left(|n|+\frac{1}{2}\right)^{\alpha}} \qquad \mbox{for certain} \quad 0 < \alpha < \frac{1}{2},$$
introducing a cut-off and making it periodic can yield functions with $r_f:\mathbb{Z} \rightarrow \left\{0,2\right\}$.
\begin{figure}[h!]
\centering
\includegraphics[width = 0.7\textwidth]{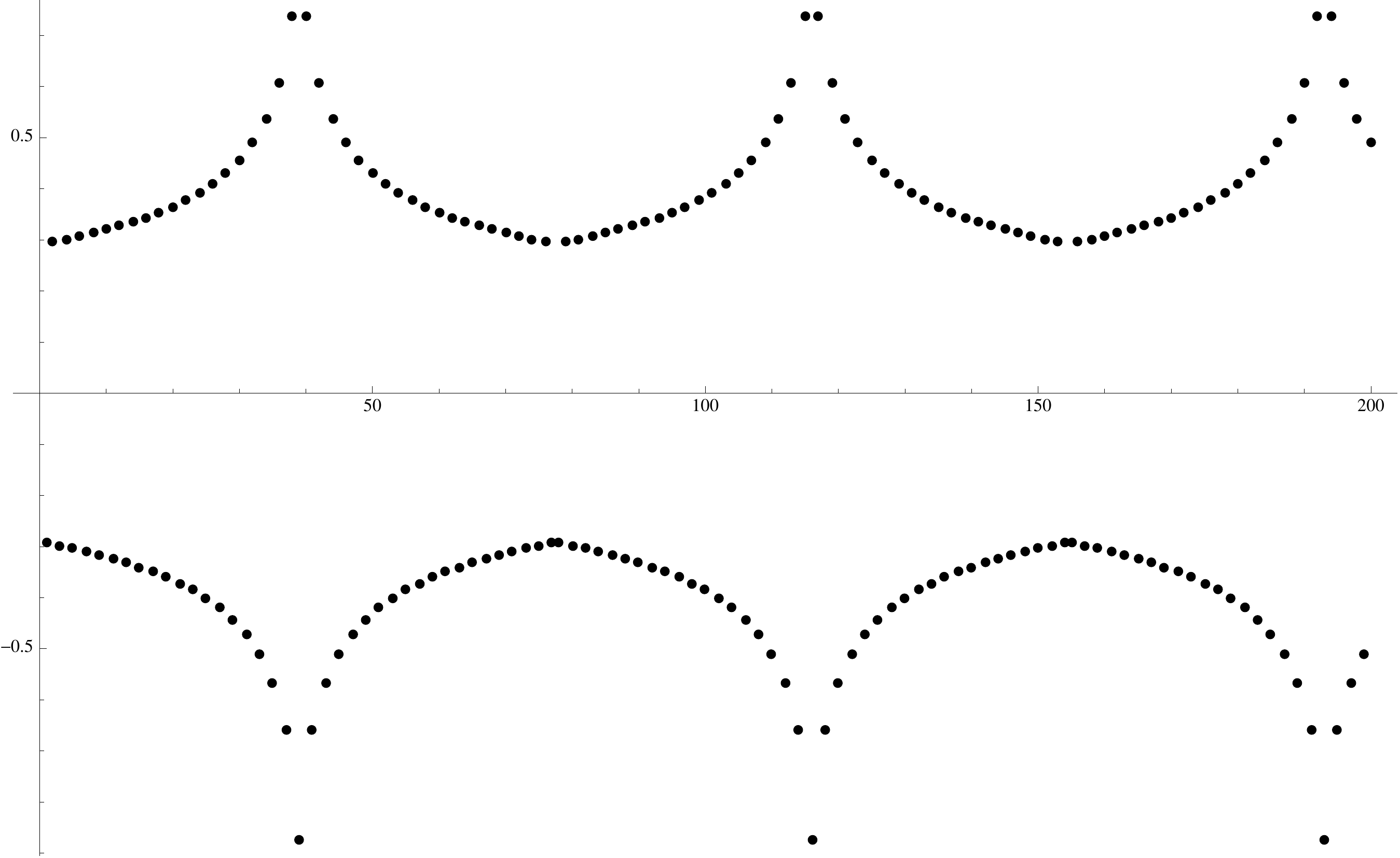}
\caption{A periodic function $f:\mathbb{Z} \rightarrow \mathbb{R}$ as above (here: $\alpha = 1/3$) with $r_f \in \left\{0,2\right\}$.}
\end{figure}

This example is structurally totally different from the sine: having oscillations at scale 2 seems crucial. We do not know whether there
are any solutions of other types and consider it to be a fascinating problem. A natural conjecture would be that if a periodic function $f:\mathbb{Z} \rightarrow \mathbb{R}$ satisfies
$$  |f(n+1) - f(n)|  \leq  \varepsilon \|f\|_{\ell^{\infty}} $$
and has a two-valued $r_f$, then
$$ \inf_{a,b,c,d \in \mathbb{R}} ~~ \sup_{n \in \mathbb{Z}} ~~\left|  f(n) - (a + b \sin{(cn + d)}) \right| \leq c(\varepsilon)\|f\|_{\ell^{\infty}}$$
for some function $c:\mathbb{R}_{+} \rightarrow \mathbb{R}_{+}$ tending to 0 as $\varepsilon$ tends to 0.

\subsection{Related work.} We are not aware of any related work in this direction. Our interest in the Hardy-Littlewood maximal function itself, however, stems
from a series of recent interesting results studying fine properties of $\mathcal{M}f$: since $\mathcal{M}f$ tries to maximize local averages, there is every reason
to believe that it should decrease total variation -- this turns out to be a surprisingly intricate problem. Motivated by a question of Kinnunen \cite{ki}, Tanaka \cite{tan} showed 
for the \textit{uncentered} Hardy-Littlewood maximal function $ \mathcal{ M^*}$ that
$$ \| (\mathcal{M}^* f)'\| \leq 2 \| f'\|_{L^1},$$
where the constant 2 was then improved to 1 by J. M. Aldaz \& J. P\'{e}rez L\'{a}zaro \cite{al}. Kurka \cite{k} has recently proven the same inequality for the centered Hardy-Littlewood maximal
function for a large universal constant. Carneiro \& Svaiter \cite{car2} give corresponding results for the maximal heat flow and the maximal Poisson operator. The discrete question on the lattice $\mathbb{Z}$
has been investigated by Bober, Carneiro, Hughes \& Pierce \cite{bob} and Carneiro \& Hughes \cite{car} in higher dimensions. The result of Kurka in the discrete setting has been proven by Temur \cite{tem}. These results are well in line of what one would expect from a maximal function, however, it is quite interesting that all of them seem quite difficult to prove; indeed, the sharp constant 1 for the centered maximal function on the real line is still merely conjectural.\\

\textit{Stokes wave.} Equations of the type appearing in our proof seems to have previously surfaced in a completely different context: in a 1987 paper on on the behaviour of the Stokes wave
of extreme form near its crest, Amick \& Fraenkel \cite{amick} encounter the equation
$$ \sqrt{3}(1+z) = \tan{\left( \frac{\pi}{2} z\right)} \qquad z \in \mathbb{C}$$
and require statements about the linear independence of the solutions of such an equation. All solutions $z_1, z_2, \dots$ with $Re(z) > -1$ are simple and the Amick-Fraenkel conjecture says that
$$ \left\{ 1, z_1, z_2, \dots \right\} \qquad \mbox{is linearly independent over the rationals}.$$
Shargorodsky \cite{shar} showed that this is implied by the Schanuel conjecture. Our proof encounters similar issues but can be unconditionally
resolved using the Lindemann-Weierstrass theorem.

\section{Proofs}
\subsection{Outline.} This section contains all the proofs. We first prove Theorem 1 and then show how it implies Theorem 2 and Theorem 3.
The proof of Theorem 1 uses an expansion into Fourier series and the fact that averaging over intervals acts diagonally on the Fourier basis. This implies that a particular Fourier series 
has to vanish identically which implies that all Fourier coefficients have to vanish identically -- this can be reduced to a system of 'diophantine' equations (over $\mathbb{N} \times \mathbb{N} \times \mathbb{R}_{+}$).
Using arguments from transcendental number theory, we can show that this system only has the trivial solution, which implies that the function has to be localized at one point in the frequency spectrum. The latter part of the argument is exploiting the arising structure in a very particular way and seems to only work in the very special case we are considering.

\subsection{Proof of Theorem 1.}
\begin{proof} Suppose that $f \in C^{\alpha}(\mathbb{R}, \mathbb{R})$ for some $\alpha > 1/2$ is periodic. Without loss of generality, we can use the symmetries of the statement to assume that the function
 has mean value 0 and the smallest period is $2\pi$ and that it can be written as
$$ f(x) = \sum_{k =1}^{\infty}{a_k \sin{k x} + b_k \cos{k x}}.$$
We assume now that
$$ \left( \partial_r  \frac{1}{2r}\int_{x-r}^{x+r}{f(z)dz} \right) \big|_{r = \gamma} = 0 \qquad \mbox{for all}~x \in \mathbb{R}.$$
The trigonometric identities
\begin{align*}
 \sin{(x+r)} - \sin{(x-r)} &= 2\sin{r}\cos{x} \\
\cos{(x+r)} - \cos{(x-r)} &= -2 \sin{r} \sin{x}
\end{align*}
 yield that
$$ \frac{1}{2r}\int_{x-r}^{x+r}{f(z)dz} = \sum_{k=1}^{\infty}{\frac{\sin{rk}}{rk}\left(a_k\sin{kx} + b_k\cos{kx}\right)}.$$
Here we invoke the classical theorem of Bernstein (see, e.g. \cite{katz}) stating that periodic functions in $f \in C^{\alpha}(\mathbb{R}, \mathbb{R})$ for some $\alpha > 1/2$ have
an absolutely convergent Fourier series. Furthermore, this allows us to interchange the sum and derivative with respect to $r$
$$0 = \partial_r  \frac{1}{2r}\int_{x-r}^{x+r}{f(z)dz}  \big|_{r = \gamma} =  \sum_{k=1}^{\infty}{\frac{\gamma k \cos{\gamma k} - \sin{\gamma k}}{\gamma^2 k}\left(a_k\sin{kx} + b_k\cos{kx}\right)}$$
because
$$ \left| \frac{\gamma k \cos{\gamma k} - \sin{\gamma k}}{\gamma^2 k} \right| \leq \frac{\gamma k + 1}{\gamma^2 k} \leq \frac{\gamma+1}{\gamma^2}$$
and therefore
$$ \sum_{k=1}^{\infty}{\left|\frac{\gamma k \cos{\gamma k} - \sin{\gamma k}}{\gamma^2 k}\left(a_k\sin{kx} + b_k\cos{kx}\right)\right|} \leq \frac{\gamma+1}{\gamma^2}\sum_{k=1}^{\infty}{\left|a_k\sin{kx} + b_k\cos{kx}\right|} < \infty.$$
The only way for a Fourier series to vanish everywhere is for all the coefficients to vanish. Note that $\gamma > 0$ and therefore
$$ \frac{\gamma k \cos{\gamma k} - \sin{\gamma k}}{\gamma^2 k} = 0 \Leftrightarrow \gamma k = \tan{\gamma k}.$$
For any fixed $k$, it is certainly possible to choose $\gamma$ in such a way that the equation is satisfied. It remains to show that no two
such equations can be satisfied at the same time. We prove this by contradiction and assume now that
$$ a_k^2 + b_k^2 > 0 \qquad \mbox{for at least two different values of}~k \in \mathbb{N}.$$
This would imply the existence of a solution $(\gamma, m, n) \in \mathbb{R} \times \mathbb{N} \times \mathbb{N}$
\begin{align*}
 \tan{\gamma m} &= \gamma m \\
 \tan{\gamma n} &= \gamma n
\end{align*}
with $\gamma > 0$ and $m \neq n$. If we could derive a contradiction from this assumption, it would imply that, independently of the value $\gamma$, 
$$ a_k^2 + b_k^2 > 0 \qquad \mbox{can hold for at most one value of}~k \in \mathbb{N}$$
from which the statement follows since then
$$ f(x) = a_k\sin{kx} + b_k\cos{kx}.$$
It is not surprising that number theory enters here: one way of rephrasing the problem is that any two elements in the set
$$  \left\{ z \in \mathbb{R}_{+}: \tan(z) = z \right\} \qquad \mbox{are linearly independent over}~\mathbb{Q}.$$
The rest of the argument can be summarized as follows: the tangent has the powerful property of sending nonzero algebraic numbers to transcendental
numbers. Any nonzero solution $\gamma \in \mathbb{R}$ of the equation $\tan{(\gamma m)} = \gamma m$ must therefore be transcendental, which means that it is never
the root of a polynomial with rational coefficients. Using multiple angle formulas for the tangent, the assumption of any nontrivial solution $\gamma$
satisfying two of these equations at the same time allows to construct an explicit polynomial for which $\gamma$ is a root -- this contradiction will conclude the proof.
We start with the cornerstone of the argument.

\begin{quote}
\textbf{Claim} (taken from \cite{mor})\textbf{.} \textit{ If $x \neq 0$ is algebraic over $\QQ$, then $\tan{x}$ is transcendental over $\QQ$.}
\begin{proof}
Suppose $\tan{x}$ is algebraic over $\mathbb{Q}$, then we would have that for some $n \in \mathbb{N}$ and some $r_k \in \mathbb{Q}$
$$ \sum_{k = 0}^{n}{r_k (\tan{x})^k} = 0.$$
We rewrite $x$ using the exponential function
$$ \tan{x} = \frac{1}{i}\frac{e^{ix} - e^{-ix}}{e^{ix} + e^{-ix}}.$$
Inserting this expression and multiplying by $(e^{ix} + e^{-ix})^n$ on both sides allows us deduce that
$$ \sum_{k = 0}^{n}{r_k \left( \frac{e^{ix} - e^{-ix}}{i}\right)^k(e^{ix} + e^{-ix})^{n-k} } = 0.$$
Expanding all brackets, we may deduce that
$$ \sum_{k = -n}^{n}{r_k^*e^{i k x} } = 0$$
for some $r_k^* \in \mathbb{Q}[[i]]$ not all of which are 0. We invoke the Lindemann-Weierstrass theorem in the formulation of Baker \cite{bak}: if 
$ b_0, b_1, \dots, b_m$ are non-zero algebraic numbers and $\beta_0, \beta_1, \dots, \beta_m$ are distinct algebraic
numbers, then
$$b_0 e^{\beta_0} + b_1 e^{\beta_1} + \dots + b_m e^{\beta_m} \neq 0$$
and this contradiction completes our proof.
\end{proof}
\end{quote}

We now prove a little statement showing that
integer multiples of fixed points $\tan{x} = x$ have a well-defined tangent.  Equivalently, we want to guarantee that if
$\tan{\gamma} = \gamma$, then $n \gamma = (m+1/2)\pi$ has no solutions $(n,m) \in \mathbb{N}^2$. 

\begin{quote}
\textbf{Claim}\textbf{.} \textit{If $\gamma > 0$ satisfies $\tan{\gamma n} = \gamma n$ for some $n \in \mathbb{N}_{>0}$, then $\gamma$ and $\pi$ are linearly independent over $\mathbb{Q}$.}
\begin{proof} Suppose that the statement fails and
$$ \gamma n =(p/q)\pi.$$
Note that, by definition,
$$ \tan\left(\tan\left( \pi p/q \right) \right) = \tan\left(\tan\left(\gamma n \right) \right) = \tan{\gamma n} =  \tan\left( \pi p/q \right),$$
It is known that $ \tan\left( \pi p/q \right)$ is an algebraic number (even the degree of the minimal polynomial is known, see \cite{ca}).
This would be an instance of the tangent mapping the nonzero algebraic number $\tan\left( \pi p/q \right)$ to an algebraic number, which
is a contradiction to the statement proven above.
\end{proof}
\end{quote}

Suppose now that $(\gamma, m, n) \in \mathbb{R} \times \mathbb{N} \times \mathbb{N}$ is a nontrivial solution of
\begin{align*}
\tan{\gamma m} &= \gamma m \\
\tan{\gamma n} &= \gamma n.
\end{align*}
Then $\gamma$ has to be transcendental: if $\gamma$ were algebraic, then $\gamma m$ would be
algebraic from which we could deduce that $\tan{\gamma m}$ is transcendental, which contradicts $\tan{\gamma m} = \gamma m$.
Now in order to derive a final contradiction exploiting the fact that $\gamma$ is transcendental 
we use an addition theorem for the tangent:
$$ \tan{((n+1) x)} = \tan{(n x + x)} = \frac{\tan{n x} + \tan{x}}{1-\tan{nx}\tan{x}}.$$
Iterating this multiple-angle formula, we have
$$ \tan{n x} = \frac{p_n(\tan(x))}{q_n(\tan(x))}$$
for two sequences of polynomials with integer coefficients satisfying the initial conditions $p_1(x) = x$ and  $q_1(x) = 1$ and
the recursion formulas
\begin{align*}
p_{n+1}(x) &= p_n(x) + x q_n(x) \\
q_{n+1}(x) &= q_n(x) - xp_n(x).
\end{align*}
We know that $(\gamma, m, n) \in \mathbb{R}_{>0} \times \mathbb{N} \times \mathbb{N}$ solves
$$ 0 = n\tan{\gamma m} - m\tan{\gamma n} = n\frac{p_m(\tan{\gamma})}{q_m(\tan{\gamma})} - m\frac{p_n(\tan{\gamma})}{q_n(\tan{\gamma})}$$
and therefore
$$ 0 = n q_n(\tan{\gamma}) p_m(\tan{\gamma}) - m q_m(\tan{\gamma}) p_n(\tan{\gamma}).$$
It is easy to see that the polynomial on the right-hand side does not vanish identically by checking that
$$ \frac{d^3}{dx^3}n \tan{x m} - m \tan{x n}\big|_{x=0} = \frac{1}{3}(nm^3 - mn^3) \neq 0.$$
This implies that $\tan{\gamma}$ satisfies a polynomial equation with integer coefficients and thus $\tan{\gamma}$ is algebraic, which is
a contradiction. 
\end{proof}

\subsection{Theorem 1 implies Theorem 2.}
\begin{proof}
  Let $f:\mathbb{R} \rightarrow \mathbb{R}$ be a periodic function of regularity $C^{\alpha}$ with $\alpha > 1/2$ such that 
$$ \left| \bigcup_{x \in \mathbb{R}}{ \left\{ r_f(x), r_{-f}(x) \right\} }  \right| \leq 2.$$
Using the symmetries of the maximal function, we may assume without loss of generality that $f$ is periodic with period $2\pi$ and
has vanishing mean value. Let us first assume that 
$$ \left| \bigcup_{x \in \mathbb{R}}{ \left\{ r_f(x), r_{-f}(x) \right\} }  \right| = 1.$$
Since $f$ is periodic, it assumes a global maximum from which it follows that if $r_f$ were to be constant, it would have to be 0
from which it follows that $f$ is constant and the statement holds. Thus we can focus on the remaining case of $r_f$ and $r_{-f}$ being two-valued
(by the same reasoning 0 always has to be one of the two values):
$$  \bigcup_{x \in \mathbb{R}}{\left\{r_f(x), r_{-f}(x)\right\}}  = \left\{0, \gamma\right\} \qquad \mbox{for some real number} ~\gamma > 0.$$
We are now in the case where $f$ is continuous, non-constant and has vanishing mean: this allows to partition $\mathbb{R}$ into three nonempty sets
\begin{align*}
I_1 &=  \left\{x \in \mathbb{R}: f(x) < 0\right\}  \\
I_2 &= \left\{x \in  \mathbb{R}: f(x) = 0\right\} \\
I_3 &=  \left\{x \in  \mathbb{R}: f(x) > 0\right\}.
\end{align*}
We will now prove that
$$ \left( \partial_r  \frac{1}{2r}\int_{x-r}^{x+r}{f(z)dz} \right) \big|_{r = \gamma} = 0 \qquad \mbox{for all}~x \in \mathbb{R}.$$
This is easy to see on $I_1$: if $x \in I_1$ and $r_{f}(x) = 0$, then the maximal function would have the value $f(x) < 0$. 
However,
by taking the maximal interval of length $2\pi$, we can at least get an average value of 0, which exceeds $f(x)$.
This implies that $r_{f}(x) = \gamma,$ which implies the statement. A similar argument works for $x \in I_3$, where 
the same reasoning implies $r_{-f}(x) = \gamma$, which implies
 $$ \left( \partial_r  \frac{1}{2r}\int_{x-r}^{x+r}{-f(z)dz} \right) \big|_{r = \gamma} = 0$$
and gives the statement after multiplication with $-1$. For $x \in I_2$ we have to argue a bit differently: suppose
$r_f(x) = 0$. Then let us consider the function $h:[0, \pi] \rightarrow \mathbb{R}$
$$h(r) = \frac{1}{2r}\int_{x-r}^{x+r}{f(z)dz}.$$
By assumption, we have that $h(0) = h(\pi) = 0$ and $h(r) \leq 0$. If $h$ vanishes identically, the derivative vanishes
everywhere and in particular also in $\gamma$. Suppose now that $h$ does not vanish identically, then it assumes a global minimum on that interval. 
By definition, this implies that $r_{-f}(x) > 0$ and thus by assumption $r_{-f}(x) = \gamma$ and this implies the statement as before -- this completes
the reduction of Theorem 1 to Theorem 2.\\
\end{proof}

\subsection{Theorem 1 implies the Theorem 3.} 
\begin{proof} Let $\gamma > 0$ be fixed and let $f \in C^1(\mathbb{R}, \mathbb{R})$ be a periodic solution of the delay differential equation
$$ f'(x+\gamma) - \frac{1}{\gamma}f(x+\gamma) = -f'(x-\gamma) - \frac{1}{\gamma}f(x-\gamma).$$
This can be rephrased as
$$ f'(x+\gamma) + f'(x-\gamma) = \frac{1}{\gamma}\left(f(x+\gamma) - f(x-\gamma)\right).$$
Integrating with respect to $x$ on both sides yields
$$ f(x+\gamma) + f(x-\gamma) = \frac{1}{\gamma}\int_{x-\gamma}^{x+\gamma}{f(z)dz} + c,$$
where $c \in \mathbb{R}$ is some undetermined constant. Since $f$ is periodic with some period $P$, we can deduce that the
average value of the left-hand side of the equation is precisely
$$ \lim_{y \rightarrow \infty}{\frac{1}{y} \int_{0}^{y}{f(x+\gamma) + f(x-\gamma) dx}} = \frac{2}{P}\int_{0}^{P}{f(x)dx}.$$
On the other hand
$$ \lim_{y \rightarrow \infty}{\frac{1}{y} \int_{0}^{y}{ \left( \frac{1}{\gamma}\int_{x-\gamma}^{x+\gamma}{f(z)dz} \right)+ c~dx}} = c + \frac{2}{P}\int_{0}^{P}{f(x)dx} $$
and thus $c=0$. Thus, multiplying the equation with $(2\gamma)^{-1}$, we get
$$ 0 = \frac{1}{2\gamma}(f(x+\gamma) + f(x-\gamma)) - \frac{1}{2\gamma^2}\int_{x-\gamma}^{x+\gamma}{f(z)dz} =  \left( \partial_r  \frac{1}{2r}\int_{x-r}^{x+r}{f(z)dz} \right) \big|_{r = \gamma}.$$
This is precisely the condition in Theorem 2 (with slightly higher regularity on $f$) and implies the result.
\end{proof}

\section{Concluding remarks}
\subsection{A conjectured stronger statement.} It seems reasonable to assume that for a periodic $C^1-$function
 $f:\mathbb{R} \rightarrow \mathbb{R}$ already
$$ \left| \bigcup_{x \in \mathbb{R}}{\left\{r_f(x)\right\}} \right| \leq 2$$
implies that $f$ has to be a trigonometric function, i.e. that it suffices to demand that the computation of $\mathcal{M}f$ is 'simple' (in the sense described above) and
not additionally that the computation of $\mathcal{M}(-f)$ be simple as well. After adding a suitable constant we can assume w.l.o.g. that $f$ has vanishing
mean and, by using the dilation symmetry, that $r_{f}(x) \in \left\{0, 1\right\}$. Then we know that $r_f(x) = 1$ whenever
$f(x) < 0$ and that therefore
$$ \partial_x \left( \partial_r \frac{1}{2r} \int_{x-r}^{x+r}{f(y) dy} \big|_{r=1} \right) = 0 \qquad \mbox{whenever}~f(x) < 0.$$
The same explicit computation as before implies that one could derive the following statement.
\begin{quote}
\textit{Conjecture.} Suppose $f:\mathbb{R} \rightarrow \mathbb{R}$
is $C^1$ and satisfies
$$  f'(x+1) - f(x+1) = -f'(x-1)-f(x-1) \qquad  \mbox{whenever}~f(x) < 0.$$
Then
$$f(x) = a+b\sin{(cx + d)} \qquad \mbox{for some}~a,b,c,d \in \mathbb{R}.$$ 
\end{quote} 

This statement would be a \textit{quite} curious strengthening of Theorem 3.

\subsection{A Poincar\'{e} inequality.}
The purpose of this short section is to note a basic observation for the uncentered maximal function (which fails for the centered maximal function): 
if the uncentered maximal intervals are all rather short, then this should imply the presence of strong oscillation in the function. We give a
very simple form of that statement. Let $f \in C^1([0,1])$ and consider the uncentered maximal function $\mathcal{M}^*$ defined via
$$  (\mathcal{M}^*f)(x) = \sup_{J \ni x}{\frac{1}{|J|}\int_{J}{f(x)dx}},$$
where $J$ ranges over all intervals $J \subset [0,1]$ containing $x$. 
We define analogously $r^*_f(x)$ as the length of the shortest interval
necessary to achieve the maximal possible value.

\begin{center}
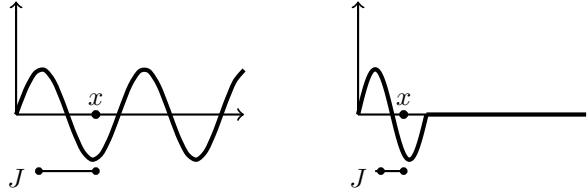
\begin{figure}[h!]
\begin{tikzpicture}[scale=1.5]
\draw [fill, ultra thick] (0.7,0) circle [radius=0.02];
\draw [fill, ultra thick] (0.7,-0.5) circle [radius=0.015];
\draw [fill, ultra thick] (0.2,-0.5) circle [radius=0.015];
  \draw [thick] (0.2,-0.5) -- (0.7,-0.5);
  \draw[->, thick] (0,0) -- (2,0);
  \draw[->, thick] (0,0) -- (0,1);
  \draw[ultra thick, scale=1,domain=0:2,smooth,variable=\x] plot ({\x},{0.4*sin(400*\x)});
 \node[below] at (0,-0.4) {$J$};
 \node[above] at (0.7,0) {$x$};

\draw [fill, ultra thick] (3.4,0) circle [radius=0.02];
\draw [fill, ultra thick] (3.4,-0.5) circle [radius=0.015];
\draw [fill, ultra thick] (3.2,-0.5) circle [radius=0.015];
  \draw [thick] (3.15,-0.5) -- (3.4,-0.5);
  \draw[ thick] (3,0) -- (5,0);
  \draw[->, thick] (3,0) -- (3,1);
  \draw[ultra thick, scale=1,domain=3:3.6,smooth,variable=\x] plot ({\x},{0.4*sin(600*\x-1800)});
  \draw[ultra thick] (3.6,0) -- (5,0);
 \node[below] at (3,-0.4) {$J$};
 \node[above] at (3.4,0) {$x$};

\end{tikzpicture}
\caption{Two examples where $r^*_f(x)$ is always small: both exhibit strong oscillation.}
\end{figure}
\end{center}

It is clear from examples that, in a loose sense, strong oscillation implies that $r^*_f(x)$ is small: it will
be optimal to choose the interval to be either a point evaluation or in such a way that one captures 
the one or two adjacent large amplitudes. An inverse result can be quantified as follows.

\begin{proposition} Assume $f \in C^1([0,1])$ has mean value $\overline{f}$. Then we have the Poincar\'{e} inequality 
$$ \int_{0}^{1}{\left| f(x) - \overline{f}  \right| dx} \leq 4\| r^*_f \|_{L^{\infty}(\mathbb{R})} \int_{0}^{1}{|f'(x)|dx}.$$
\end{proposition}
If we were to replace $4\|r^*_f\|_{L^{\infty}}$ by the constant $1/2$, this would be the classical $L^1-$Poincar\'{e} inequality on $[0,1]$. 
Put differently, for a function $f:[0,1] \rightarrow \mathbb{R}$ with vanishing mean, we have 
an uncertainty relation between the total variation and $ \|r^*_f\|_{L^{\infty}}$
$$ \mbox{var}(f) \|r^*_f\|_{L^{\infty}} \geq \frac{1}{4}\int_{0}^{1}{|f(x)|dx}.$$
It is easy to see that the statement has the sharp scaling: consider 
$$ f(x) = \sin{N \pi x}  \quad \mbox{where} \qquad \|r^*_f \|_{L^{\infty}} \sim N^{-1} ~ \mbox{and} ~\| f'\|_{L^1} \sim N.$$
Another example is given by taking a positive bump function $\phi \in C^{\infty}_{c}(0,1)$ and consider 
the rescaled function
$$ f(x) = a \phi'(b x) \quad \mbox{where} \qquad \int_{0}^{1}{\left| f(x) - \overline{f}  \right| dx} \sim a b^{-1}, \|r^* \|_{L^{\infty}} \sim b^{-1} ~ \mbox{and} ~\| f' \|_{L^1} \sim a.$$
We note that the result is not true for the centered maximal function: the function $f(x) = 1 - (x-0.5)^2$ (or, more generally,
any strictly concave function) satisfies $r_f \equiv 0$.

\begin{proof} We may suppose w.l.o.g. that $f$ has vanishing mean $\overline{f} = 0$. Now we write
 $$\int_{0}^{1}{\left| f(x)  \right| dx} = 2\int_{0}^{1}{\chi_{f < 0}(x)|f(x)| dx}$$
and estimate the number on the right. It is easy to see that the set 
$$ \left\{x \in [0,1]: f(x) < 0\right\}$$
cannot contain an interval of length larger than $2\|r^{*}_f\|_{L^{\infty}}$ because, since $\overline{f} = 0$, the maximal function is always nonnegative
(one can always simply choose the entire interval). Clearly, for any $g \in C^1$, 
$$ \int_{s}^{t}{|g(x)|dx} \leq (t-s)\int_{s}^{t}{|g'(x)| dx} \qquad \mbox{if}~g(s)=0.$$
We can now use that inequality on every connected component $ \left\{x \in [0,1]: f(x) < 0\right\}$: by the reasoning above, we will always have $t-s \leq  2\|r^{*}_f\|_{L^{\infty}}$
and therefore
$$\int_{0}^{1}{\chi_{f < 0}|f| dx} \leq 2\|r^{*}_f\|_{L^{\infty}} \int_{0}^{1}{| f'(x)|dx},$$
which implies the statement.
\end{proof}
It could be quite interesting to 
understand under which conditions and to which extent such improved Poincar\'{e} inequalities are
 true in higher dimensions and how they depend on the maximal function involved.\\

\textbf{Acknowledgement.} I am grateful for discussions with Raphy Coifman and indebted to Lillian Pierce for a number helpful remarks, which
greatly improved the manuscript. My interest in non-standard aspects of the maximal function traces back to a series of very enjoyable conversations 
with Emanuel Carneiro at the Oberwolfach Workshop 1340 -- I am grateful to both him and the organizers of the workshop.

\end{document}